\pdfoutput=1
\documentclass[a4paper,reqno,oneside,11pt]{article}
\usepackage[utf8]{inputenc}
\usepackage{graphicx}
\usepackage{float}
\usepackage{amsmath,amscd}
\usepackage{amsfonts}

\usepackage{amssymb}
\usepackage{amsthm}
\usepackage{epsfig}
\usepackage{mathtools}
\usepackage[english]{babel}
\usepackage{textcomp}
\usepackage[T1]{fontenc}
\usepackage{dsfont}

\usepackage{tikz-cd}
\usepackage{tikz}
\usepackage{pgfplots}
\usetikzlibrary{decorations.markings}
\usetikzlibrary{shapes.geometric, arrows}
\usepackage{tikzscale}
\usepackage{filecontents} 
\usepackage{wrapfig}
\usepackage{tcolorbox}

\makeatletter
\newcommand{\mybox}[1]{  \setbox0=\hbox{#1}  \setlength{\@tempdima}{\dimexpr\wd0+13pt}  \begin{tcolorbox}[colframe=gray1,boxrule=0.5pt,arc=4pt,
      left=6pt,right=6pt,top=6pt,bottom=6pt,boxsep=0pt,width=\@tempdima]
    #1
  \end{tcolorbox}
}
\makeatother

\usetikzlibrary{arrows}

\tikzset{middlearrow/.style={
        decoration={markings,
            mark= at position 0.4 with {\arrow{#1}} ,
        },
        postaction={decorate}
    }
}

\tikzstyle{arrow} = [ very thick,->,>=stealth]
\usepackage{caption}

\definecolor{brown}{RGB}{150,100,0}

\definecolor{MISred}{rgb}{0.533 0.085 0.22}
\definecolor{MISblue}{rgb}{0.18 0.353 0.52}
\usepackage{hyperref}
\usepackage{subfig}
\usepackage{calc}
\usepackage{picinpar}
\usepackage[all,cmtip]{xy}
\usepackage{epstopdf}
\usepackage{newfile}
\usepackage[a4paper]{geometry}

\usepackage{bbm}
\usepackage{subfig}
\usepackage{float}
\usepackage{enumitem}
\setlist[enumerate,2]{ref=\theenumi(\alph*)}

\newlist{axioms}{enumerate}{1}
\setlist[axioms,1]{
  label={Axiom~\arabic*},
  leftmargin=*,
  align=left
  }

\usepackage{calc}

\usepackage{fix-cm}

\usepackage{color}
\definecolor{gray}{rgb}{0.9,0.9,0.9}
\definecolor{gray1}{rgb}{0.7,0.7,0.7}
\definecolor{gray2}{rgb}{0.8,0.8,0.8}
\definecolor{magenta}{rgb}{1.0, 0.0, 1.0}

\hypersetup{
  colorlinks   = true,   urlcolor     = blue,   linktoc =  page,
  linkcolor    = red,   citecolor   = green }
\definecolor{malachite}{rgb}{0.04, 0.85, 0.32}

\usepackage{etoolbox}

\usepackage{booktabs,tabularx}
\usepackage{microtype}
\allowdisplaybreaks

 \newoutputstream{pages}

\openoutputfile{pages.num}{pages} 
\addtostream{pages}{\arabic{page}}
\closeoutputstream{pages}

\theoremstyle{plain}
\newtheorem{Theorem}{\Large T\normalsize heorem}[section]
\newtheorem{Proposition}[Theorem]{\Large P\normalsize roposition}
\newtheorem{Lemma}[Theorem]{\Large L\normalsize emma}

\newtheorem{Definition}{\Large D\normalsize efinition}[section]
\theoremstyle{definition}
\newtheorem{Example}{Example}[section]
\newtheorem{Remark}{Remark}[section]

\newcommand{\bel}[1]{\begin{equation}\label{#1}}

\newcommand{\be}{\begin{equation}} 
\newcommand{\qe}{\end{equation}}
\newcommand{\Se}{\mathbb{S}}     

\newcommand{\ba}{\begin{eqnarray}}
\newcommand{\ea}{\end{eqnarray}} 
\newcommand{\rf}[1]{(\ref{#1})}
\newcommand{\bi}{\bibitem}

\DeclareMathOperator{\im}{im}

\DeclareMathOperator{\Proj}{Proj}
\DeclareMathOperator{\ind}{ind}

\newcommand{\si}{\sigma}

\newcommand{\K}{\mathbb{K}}
\newcommand{\R}{\mathbb{R}}
\newcommand{\Z}{\mathbb{Z}}
\newcommand{\bo}{\partial}

\title{Discrete Morse-Bott theory for CW complexes}
\author{Sylvia Yaptieu  \thanks{Max Planck Institute for Mathematics in the Sciences, Inselstrasse 22, 04103 Leipzig, Germany.}}

\begin{document}

\maketitle

\begin{abstract}
We derive a discrete analogue of Morse-Bott theory on CW complexes and  use this discrete Morse-Bott function to 
 do some  Conley theory analysis. 
      It turns out that our discrete Morse-Bott theory is indeed a generalization of Forman's discrete Morse theory in \cite{For1}. 
       
 \end{abstract}
2010 Mathematics Subject Classification 37F20, 37B30, 55N99.

\section{Introduction}

Morse theory is a very important tool for the study of the topology of differentiable manifolds. It was first introduced by
Morse in 1925, in \cite{MMo}.
It recovers the homology of the manifold from the critical points of a Morse function and the 
relations between them. 
The Morse inequalities are inequalities between the Betti numbers  of the manifold and the numbers of critical
 points of fixed indices of the function. 
To get the homologies, one attaches a $k$-dimensional cell for each critical point of index $k$, and  gluing relations
between those cells then 
yield the homological boundary operator. 

Morse-Bott theory, introduced by Bott in \cite{RB2}, being a generalization of Morse theory, was  developed to treat the cases 
where instead of having critical points, one has critical submanifolds. To each critical submanifold, one associates  a certain index 
that is determined by looking at the Hessian restricted to the normal part of the submanifold, since on the tangential part it vanishes. 
The Poincar\'e polynomial (and hence the Euler number) of the manifold is then obtained from the Morse-Bott inequalities in the sense that it 
is expressed in terms 
of the Poincar\'e polynomials of the critical submanifolds, taking their respective indices into account.

Conley theory on the other hand, introduced by Conley in \cite{CC}, and being more dynamics related, focuses  on the study of the topological invariants 
of a given manifold. Using the (negative) gradient flow lines generated by some (Morse) function, 
one  obtains for each isolated invariant set  
its isolating neighborhood and exit set for the flow, and these two constitute an index pair. One obtains the Poincar\'e polynomial of 
the manifold by summing up those of the index pairs up to some correction term.
 In particular, 
 the Euler number of the manifold is then obtained by summing for all isolating invariant sets the alternating sums of the dimensions 
of the homology groups of  the index pairs.

Forman in 1998, in \cite{For1}, developed a discrete version of Morse theory for   CW complexes. 
A CW complex, introduced by J. C. Whitehead in \cite{JCH}, is a decomposition of a space into cells each of which is homeomorphic to an open disc.  
 The dimension of the disc specifies the dimension of the cell. The CW construction  uses a specific gluing procedure via the characteristic maps. 
 This space  is endowed with the weak 
topology and satisfies the condition that the closure of each cell  intersects only a finite number of cells. We write $\si^{(k)}$ to emphasize that $\si$ is a cell 
of dimension  $k$. 
The topological boundary elements of a cell are called its faces. If a  cell $\si^{(k)}$ is a face of another cell $\tau$, 
we write $\sigma < \tau$ to indicate that 
 $\,\dim \si=\dim \tau-1$, in which case $\si$ is called a facet of $\tau$.  We say $\si^{(k)}<\tau$ is  a regular facet of  $\tau$ 
 if, for $\varphi_{\tau}$,
the characteristic map of $\tau$,  we have: the map 
 $\varphi_{\tau}\colon\varphi^{-1}_{\tau}(\sigma^{(k)})\,\rightarrow \, \sigma^{(k)}$ is a homeomorphism and 
$\overline{\varphi^{-1}_{\tau}(\sigma^{(k)}) }$ is a closed $k$-ball.

We denote the cardinality of a set $A$ by $\sharp A$.

   A discrete Morse function, according to Forman, is a real-valued function defined on the set of cells such that it locally increases in dimension, 
   except possibly in one direction. More formally, 
 Forman's definition of a discrete Morse function $f$ on a CW complex requires that for all cells 
$\sigma^{(k)}$,  
\ba\nonumber\begin{cases}
\text{for all}\,\, \tau \,\,\text{s.t.}\,\,\si\,\,\,\text{is an irregular face of}\,\, \tau, \,\,\, f(\si)< f(\tau);\\ 
  Un(\si):=\sharp \{\tau^{(k+1)}\,|\, \si \,\,\text{is a regular facet of}\,\, \tau \,\text{and}\, f(\tau) \leq f(\sigma)\} \leq 1; 
\end{cases}
\ea
and 
\ba\nonumber\begin{cases}
\text{for all}\,\, \nu \,\,\text{s.t.}\,\,\nu\,\,\,\text{is an irregular face of}\,\, \si, \,\,\, f(\nu)< f(\si);\\ 
  Dn(\si):=\sharp \{\nu^{(k-1)}\,|\, \nu \,\,\text{is a regular facet of}\,\, \si \,\text{and}\, f(\nu) \geq f(\sigma)\} \leq 1.
\end{cases}
\ea
When $\si$ is not a regular facet of some other cell, it is automatically critical; in the regular case,
$\sigma$ is critical if both  $Dn(\si)$ and $Un(\si)$ are $0$. In fact, one easily sees that at most one of 
them can be $1$; the other then has to be $0$.

A CW complex is regular if all the faces are regular. As examples one has polyhedral, cubical and simplicial complexes.

A pair $\{\si,\tau\}$ with $\si<\tau$ and $f(\si)\geq f(\tau)$ is called a noncritical pair.
If we draw an arrow from $\si$ to $\tau$ whenever $\si<\tau$ but $f(\si)\geq f(\tau)$, then we get a vector field  associated to this function, 
and  each noncritical cell 
has precisely one arrow which is either incoming or outgoing. Therefore, for the Euler number,  we only need to count the critical cells 
with appropriate signs according to their dimensions, since the noncritical cells cancel in pairs. 

We recall that, see \cite{For3}, a  combinatorial vector field on a CW complex $\K$ is a map
$\,\, V\colon \K \rightarrow \K \,\cup \, \{ 0\}\,\, $ that satisfies: if the image of a cell is non zero, 
then the dimension of the image is the dimension of the cell plus one; 
the image cells have zero images; for each cell $\si$, either $V(\si)=0$ or $\si$ is a regular face of $V(\si)$; 
and the pre-image of a given cell contains at most one element.
One thinks of the function $V$ as assigning an arrow from $\si$ to $\tau$ whenever $V(\si)=\tau$. 
In this way, each cell can either have one incoming or outgoing arrow but not both.
We write $\si \rightarrow \tau$ whenever there is an arrow from $\si$ to $\tau$.

 It is an important consequence of Forman's definition that the vector field admits no closed orbits, where by a closed orbit we mean  
 a path of the form 
$$ \sigma_0 \rightarrow \tau_0 > \sigma_1 \rightarrow \tau_1 > \dots \si_m \rightarrow \tau_m > \sigma_0.$$
Since a combinatorial vector field can also admit closed orbits, 
 the vector field constructed from a discrete Morse function is a combinatorial vector field that has no closed orbits. Conversely, one can always construct 
a discrete Morse function from a combinatorial vector field that admits no closed orbits.

Forman answered the question of a  discrete analogue of Conley theory for CW complexes, see \cite{For3}.  He 
first  uses a combinatorial vector field, and as isolated invariant sets he considers the rest points (which are
 the critical cells) and the closed orbits. The isolating neighborhoods here are the unions of all the cells 
in the isolated invariant sets together with the ones  in their boundaries; the exit set is just the collections of cells in 
the isolating neighborhood that are not 
 in the isolated invariant sets.

 Our  goal is to construct a discrete analogue of Morse-Bott theory. This involves generalizing discrete Morse functions to include the case where a discrete function 
can have 
 larger critical collections of cells, instead of only simple critical cells. 
    We define a discrete Morse-Bott function 
 by requiring  some conditions on specific collections of cells, where the function assumes the same value for all the cells in the collection.
  Using the reduced collections, excluding the noncritical pairs, we obtain 
some discrete Morse-Bott inequalities.  That is, the Poincar\'e polynomial  
of the CW complex is expressed in terms of those of the reduced collections. 
The vector field originating from this discrete Morse-Bott function is such that  
inside each collection a cell can have possibly more than one  incoming and/or outgoing arrow, but between the collections 
there are no closed orbits.
 We also do some Conley theory analysis by using the reduced collections (excluding the noncritical pairs) 
 as the isolated invariant sets, and define their respective isolating neighborhoods 
    and exit sets. These two constitute the index pairs.
    The Poincar\'e  polynomial of the CW complex is then expressed in terms of those of the index pairs.
    
    It should be noted that, when considering our discrete Morse-Bott function, 
the extracted discrete  vector field 
is not a combinatorial vector field, thus Forman's discrete Conley theory 
and our approach are complementary.

\section{Discrete Morse-Bott Theory}\label{sec5.3}

We recall, see \cite{RB2}, \cite{BH}, that in the smooth setting, a Morse-Bott function 
is one that admits critical submanifolds instead of only isolated nondegenerate critical points as
it is required for a Morse function (see \cite{Sch}, \cite{JM} and \cite{MMo}). The Poincar\'e polynomial of the manifold is retrieved by adding those of the 
critical submanifolds taking their respective indices into account, up to some correction term.\\

As mentioned earlier, an analogue of Morse-Bott theory in the discrete case might be of importance, 
since  extending a function defined on the set of vertices might not
always result in a discrete Morse function everywhere on the cell complex. The idea is, we consider a function that is   
discrete  Morse 
on a cell complex except possibly at some maximal  collection of cells, where every cell in the collection has the same value. 
We derive the analogue of the smooth Morse-Bott inequalities, that is,  
 inequalities involving the Poincar\'e polynomial (and hence the Betti numbers) of the cell complex and those of the 
 reduced collections.
First we shall need an analogue of the Morse-Bott function in the discrete setting.
\begin{Definition}[Collection]
 Let $f$ be a discrete function on some CW complex. A collection $C$ for $f$ is a maximal set of cells such that:
 \begin{itemize}
  \item [$(i)$] all the cells in $C$ have the same function value,     \item[$(ii)$] $\bigcup_{\si \in C} \bar{\si}$ is (path-) connected.
 \end{itemize}

\end{Definition}
We recall that a pair $\,\{ \si\,,\,\tau\}\,$ is said to be  noncritical if $\,\si< \tau\,$ and $\,f(\si)\geq f(\tau)$.

\begin{Definition}[Discrete Morse-Bott function]\label{dmbtfnc} Let $f$ be a function defined on a CW complex $\,\K$ and let $\,C^1,\cdots,C^l$ 
be the collections for $f$. We say $f$ is  discrete Morse-Bott if, for each $i$ and 
 for all $\,\, \si^{(k)}\in C^i,\,\,$

if $\si$ is an irregular face of some  $\tau$, $f(\si)<f(\tau)$, else 
 \ba\label{coni}
U^c_n(\si)&:=&\sharp\{ \tau^{(k+1)} \notin C^i\,| \,\,\si\,\,\text{a regular facet of}\,\,\tau, \,\, f(\tau)< f(\si)\} \leq 1;
\ea

if $\nu<\si$ is an irregular face of $\si$, $f(\nu)<f(\si)$, else 
\ba \label{conii}
D^c_n(\si)&:=& \sharp\{ \nu^{(k-1)}\notin C^i\,| \,\nu\,\, \text{a regular facet of}\,\,\si,\,\, f(\nu)> f(\si)\} \leq 1.
 \ea

 \end{Definition}
 \begin{Remark}
\begin{itemize}
 \item[$(i)$] If for all $i$, $C^i=\{\si_i\}$ and  $f$ is a discrete Morse-Bott function, then $f$ is also a discrete  Morse function.
\item[$(ii)$] Every discrete Morse function is a discrete Morse-Bott function in which each collection has at most two elements.
\end{itemize}
 \end{Remark}
\begin{Remark}
 From the definition above, we get that in any collection there cannot be 
 cells $\si$ and $\tau$ with $\si$ an irregular face of $\tau$. This is useful in the situation where, if a collection reduces to only two cells of adjacent dimension, then 
 one should not be an irregular facet of the other. This is important because in discrete Morse theory, 
 the noncritical pairs always have no contribution.
\end{Remark}
It can be shown that either $U^c_n(\si)=1$ or $D^c_n(\si)=1$, but not both. The argument used is the same as in the discrete Morse case, by considering a discrete 
 Morse-Bott function for which all the collections are singletons.\\

\begin{figure}[htbp]
\begin{center}

 \includegraphics[width=0.9\textwidth]{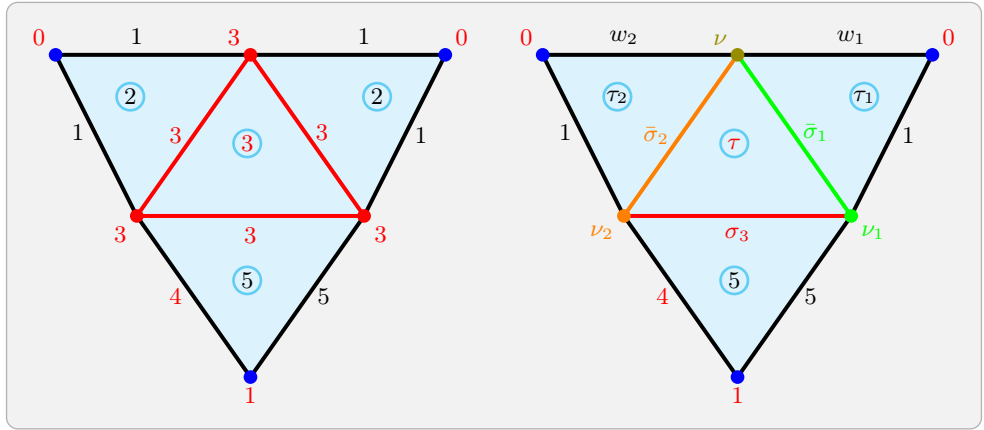}
\caption{A discrete function that is not Morse-Bott.}\label{fig:mbt}
\end{center}
\end{figure}

Figure \ref{fig:mbt} shows an example of a discrete function  that is not a discrete Morse-Bott function. 
Take $C=\{\nu,\nu_2,\nu_1,\si_1,\si_2,\si_3,\tau \}$, then the vertex $\nu$ violates 
 \rf{conii}, indeed there are  $w_1$ and $w_2$ not in $C$, with $w_1>\nu<w_2$ but 
$f(w_2)<f(\nu)>f(w_1)$.

\begin{Definition}
 Let $f$ be a discrete Morse-Bott function and $C$  a collection. 
We say that $\si \in C$ is \textbf{upward noncritical with respect to $C$} if there exists a cell  
$w(\si) \notin C$, $\,\,w(\si)>\si$  s.t. $\,\,f(\si)> f(w(\si))$.
\end{Definition}
\begin{Definition}
 Let $f$ be a discrete Morse-Bott function and $C$  a collection. 
We say that $\si \in C$ is \textbf{downward noncritical with respect to $C$} if there exists a cell  
$w(\si) \notin C$, $\,\,w(\si)<\si$  s.t. $\,\,f(\si)< f(w(\si))$.
\end{Definition}

\begin{Remark}
 If $C$ is a singleton in the two definitions above, then we have the usual (upward and downward) noncriticalities in  Forman's framework.
\end{Remark}

Observe that a cell cannot be at the same time upward and downward noncritical with respect to the same collection. 

The following is the analogue of a critical submanifold.
\begin{Definition}[Reduced collection]\label{redcol}
For each collection $C$, we define the reduced collection $C^{red}$ by taking out of $C$ all the cells that are upward or downward noncritical 
with respect to $C$. 
\end{Definition}
\begin{Remark}
 \begin{itemize}
  \item [$1)$] If $C=\{\si\}$ with $\si$ neither upward nor downward noncritical, then $\si$ is critical in the discrete Morse sense for the function $f$ and $C=C^{red}$.
\item[$2)$] If $C=\{\si\}$ where $\si$ is either upward or downward noncritical, then $C^{red}=\emptyset$.

\item[$3)$] Observe that if $\,C^{red}$ consists of only one element, then this element is not necessarily critical in the usual sense for the function $f$. 
To see this just consider 
 a given function and a  $C=\{\si_1, \nu, \si_2\},$ with $\si_1>\nu<\si_2$ but both $\si_1$ and $\si_2$ 
 being  downward noncritical. 
Then $C^{red}=\{\nu\},$ but $\nu$ is not critical in the usual sense for this function since it has the 
same value with $\si_1$ and $\si_2$.
 \end{itemize}

\end{Remark}

The following definition helps to distinguish between a cell that we will call critical and a reduced collection that is a singleton.

\begin{Definition}[Critical cell]
 A cell $\si$ is said to be critical if 
 $$\sharp \{ \tau > \si | f(\tau)\leq f(\si)\}=0\quad \text{and}\quad \sharp \{ \nu<\si | f(\nu)\geq f(\si)\}=0.$$
\end{Definition}

The definition above tells us that  if $C=C^{red}=\{\si\}$ then $\si$ is critical.

The following lemma establishes the fact that the faces of upward noncritical cells w.r.t. a collection are also upward noncritical.

\begin{Lemma}\label{Lem01}
 Let $C$ be a collection for a discrete Morse-Bott function. If $\si \in C$ is  
 upward noncritical with respect to $C$, then every face of $\si$ in $C$ is also upward noncritical w.r.t. $C$.\\
 In particular, if $C$ is a subcomplex and $\si \in C$ is  
 upward noncritical with respect to $C$, then so is $\bar{\si}$.
\end{Lemma}
\begin{proof}  Let $\,\,\nu<\si,$ and $\nu\in\,C$. The fact that $\,\,\si$ is upward noncritical w.r.t. $C$ implies there is 
 an $\,w(\si) \notin C,\,\,$ $w(\si)\in C'$ with  
 $\,w(\si) > \si,\,$
 s.t. $\,f(\si)> f(w(\si))$. Since there are no irregular faces in a collection we get that $\nu$ is a regular face of $\si$, 
 and the fact that there are no closed orbits  tells us that $\si$ is a regular face 
 of $w(\si)$. We have 
 $\nu<\si<w(\si)$, that is $\nu$ and $w(\si)$ are incident. 
 Let $\tau \notin C$ be such that $w(\si)> \tau >\nu<\si,\,\,$ such a $\tau$ exists by the incidence property.
 Then either $\tau \in C'\,\,$ or $\,\,\tau \notin C'$: if $\,\tau \in C'$ then $f(\tau)=f(w(\si))$;  
 if $\,\tau \notin C'$, the definition of a discrete Morse-Bott function yields $f(w(\si))> f(\tau)$. Thus, 
$f(\nu)= f(\si)> f(w(\si))\geq f(\tau)$, therefore  $\nu$ is also upward noncritical w.r.t. $C$, hence  
all the faces of $\si$ in $C\,\,$
  are also upward noncritical with respect to $C$.
\end{proof}

The most important fact about our definition of a discrete Morse-Bott function is the following. 
\begin{Lemma}\label{Lem02}
 If a collection $C$ is such that $C\neq C^{red}$, then there is a one-to-one correspondence between the upward (resp. downward) noncritical 
cells  $\si$ w.r.t. $C$ and the cells $w(\si) \notin C,\,\,$ $w(\si)>\si$ with $f(\si)> f(w(\si))$ 
(resp. $w(\si)<\si$ with $f(\si)< f(w(\si))$).
\end{Lemma}
\begin{proof}
We already know from Lemma \ref{Lem01} that if $\si \in C$ is  
 upward noncritical with respect to $C$, then every (regular) face of $\si$ in $C$ is also upward noncritical w.r.t. $C$. 
 If $C$ is a subcomplex,   the interesting case would be to consider a situation where 
 $\si_1, \,\si_2 \in C$ are upward noncritical w.r.t. $C$ 
 satisfying $\bar{\si}_1 \nsubseteq \bar{\si}_2$ and $\bar{\si}_2 \nsubseteq \bar{\si}_1$. 
     Then, 
 if we take a $\nu \in C\,$ such that $\,\si_1>\nu<\si_2$,  
we  have the existence of other cells $\,\,w(\si_i)\notin C\,\,$ for $\,i=1,2\,\,$ s.t. $\,\,w(\si_i)>\si_i$. 
Thus $\nu$ is a regular face of the $\si_i$'$s$ and 
this also implies (by the incidence property), there exist 
$\tau_i \notin C\,\,$ for $i=1,2\,$ with  $\,\,\nu<\tau_i<w(\si_i),$  
 and 
the following holds:
$$f(\nu)=f(\si_1)> f(w(\si_1))\geq f(\tau_1) \quad \text{and} \quad f(\nu)=f(\si_2)> f(w(\si_2)) \geq f(\tau_2).$$

Thus, if $\tau_1 \neq \tau_2$, $f$ at $\nu$ will have a value greater than that of more than one cell in 
the complement of $C$, and 
this will contradict the definition of a discrete Morse-Bott function. 
If however $\tau_1=\tau_2$, then there is a unique $w(\nu)=\tau \notin C$ for which $f(\nu)>f(w(\nu))$.

If $\si \in C$ is downward noncritical, then it follows directly from the definition of $f$.
\end{proof}
Figure \ref{fig:mbt} also shows  that the result in Lemma \ref{Lem02} does not hold whenever we have a function that is 
not discrete Morse-Bott.

The notion of the index here however is ambiguous, but it will be shown in subsequent examples 
 that the Euler number will always be counted with positive sign. See for instance Figure \ref{fig:expl4}, Figure \ref{fig:expl2}, and Figure \ref{fig:expl3}; 
  each of them involves a simplicial complex $\K$ with
negative, zero and positive Euler number, and also reduced collections with negative, 
zero and positive Euler numbers. Recall that the contribution of a cell of index $k$, in the 
computation of the Euler number of a CW complex, is given by $(-1)^k$.

\begin{figure}[htbp]
\begin{center}

 \includegraphics[width=0.9\textwidth]{figure/expl4.tikz}
\caption{A discrete Morse-Bott function on $\K$ with $\chi(\K)=-1$.}\label{fig:expl4}
\end{center}
\end{figure}

\begin{figure}[htbp]
\begin{center}

 \includegraphics[width=0.9\textwidth]{figure/expl2.tikz}
\caption{A discrete Morse-Bott function on $\K$ with $\chi(\K)=0$.}\label{fig:expl2}
\end{center}
\end{figure}

\begin{figure}[htbp]
\begin{center}

 \includegraphics[width=0.9\textwidth]{figure/expl3.tikz}
\caption{A discrete Morse-Bott function on $\K$ with $\chi(\K)=1$.}\label{fig:expl3}
\end{center}
\end{figure}

Given that the reduced collections are not necessary subcomplexes, we need to make precise how their Betti numbers are obtained. First we prove the following:\\
\begin{Lemma}\label{SbCb}
 If $\si\in \overline{C^{red}}\setminus C^{red}$ and $f(\si)=f(\tau)\,$ for $\,\tau\in C^{red}$, then $\,\si$ cannot be downward noncritical.
\end{Lemma}
\begin{proof}
 If $\si\in \overline{C^{red}}\setminus C^{red}, \,\,$ then $\,\, \si\notin  C^{red}\,$ but $\,\si$ is a face of an element in $C^{red}$, 
say $\si<\tau$ for some $\tau\in C^{red}$. \\
If $\si$ is downward noncritical, then there is a $\nu^{\star}<\si$ s.t. $f(\nu^{\star})>f(\si)$ and all 
  the other $(\nu\notin C^{red})$ $\nu<\si$ are such that $f(\si)>f(\nu)$. Note that $\nu^{\star}$ and $\tau$ are incident. Thus there is a 
$\,\,\si^{\star}>\tau\,\,$ s.t.$\,\,\nu^{\star}<\si^{\star} \neq \si$.\\
If $\si^{\star}\in C^{red}\,$ we have a contradiction since $\,f(\nu^{\star})>f(\si)=f(\si^{\star})\,$ that is, $\,f(\nu^{\star})\,$ 
will be greater than $f(\si)$ and $f(\si^{\star})$.
Thus $\,\si^{\star}\notin C^{red},$ and  $f(\si^{\star})>f(\nu^{\star})$ because $\nu^{\star}$ is already upward noncritical with $\si$. 
This automatically means that $f(\si^{\star})>f(\nu^{\star})>f(\si)=f(\tau),$ which also contradicts the fact that $\tau \in C^{red}$.
\end{proof}

 \begin{Lemma}\label{SbC}
  For any reduced collection $C^{red}$, $\,\,\overline{C^{red}}\setminus C^{red}$ is a subcomplex.
 \end{Lemma}
\begin{proof}
 
Let $\si\in \overline{C^{red}}\setminus C^{red}$, and $\nu<\si$, 
then by definition of $\overline{C^{red}},$ $\,\nu$ is also in $\overline{C^{red}}$. 
Thus, to show that $\nu\in \overline{C^{red}}\setminus C^{red}$,  we only need to show that $\nu\notin C^{red}$. 

For $\si\in \overline{C^{red}}\setminus C^{red}, \,\,$  $\si$ is a face of an element in $C^{red}$, 
say $\si<\tau$ for some $\tau\in C^{red}$. 
 The fact that $\tau\in C^{red}$ and $\si<\tau$ implies that $f(\si)\leq f(\tau)$.\\
 If $f(\si)=f(\tau)$, this means that $\si$ is either downward or upward noncritical w.r.t. $C$.
>From Lemma \ref{SbCb}, such a $\si$ cannot be downward noncritical.
 If $\si$ is upward noncritical, then since $\si$ cannot be downward noncritical, we have $f(\nu)\leq f(\si)$. 
If $f(\nu)=f(\si)$, then $\nu$ has to be in $C$, and by Lemma \ref{Lem01}, $\nu$ 
also has to be upward noncritical w.r.t. $C$, which implies that $\nu\notin C^{red}$. 
If $f(\nu)<f(\si)$ then $\nu\notin C^{red}$ by definition of $C^{red}$.\\
If 
$f(\si)< f(\tau)$  then  either $f(\nu)\leq f(\si)<f(\tau)$, $\,f(\si)< f(\nu)<f(\tau) $, or $\,f(\nu)=f(\tau)>f(\si)$. 
 In any case $\nu$ cannot  be in $C^{red}$.
\end{proof}
\begin{Remark}
 For a given collection $C$, $\,\overline{C^{red}}$ is not always equal to $C$. Figure \ref{fig:mbc1} illustrates this: $C$ is the 
 collection of the $2$-cell with value $3$, the two red edges, the red vertex and the two green vertices. After removing the upward noncritical 
 vertices w.r.t. $C$, that is the two green vertices, and the downward noncritical $2$-cell w.r.t. $C$, we end up with the two red edges and the red vertex. 
 Thus, $\overline{C^{red}}$ consists of the two red edges, the red vertex and the two green vertices, which is different from $C$.
\end{Remark}

Let $[\tau:\si]$ denote the incidence number of $\tau$ and $\si$, which is the number of times that $\tau$ is wrapped around $\si$, taking the induced orientation 
from $\tau$ onto $\si$ into account. \\

Let $C_k(C^{red};\Z)$ be the free $\Z$-module generated by the oriented cells of $C^{red}$ of dimension $k$. One can also take $C_k(C^{red};\Z_2)$.
\begin{Definition}\label{Bored}
 The boundary operator $\,\bo^{red}_k\colon C_k(C^{red},\Z)\rightarrow C_{k-1}(C^{red},\Z)$ is given by: 
 $$\bo^{red}_k\tau^{(k)}=\sum_{\si \in C^{red},\si<\tau}[\tau:\si] \si^{(k-1)}.$$
\end{Definition}
\begin{Proposition}
 $$\bo^{red}_{k-1} \circ \bo^{red}_k=0.$$
 \end{Proposition}
 \begin{proof}
  If $C^{red}$ is a subcomplex then it follows from the classical theory. \\
  If $C^{red}$ is not a subcomplex, the result will follow if we show that the above boundary operator is a relative boundary operator, that is,  
a boundary operator for relative homology.
   Indeed\\
  
  $\overline{C^{red}}:=\bigcup_{\si \in C^{red}}\, \bar{\si}$ is a subcomplex by definition.\\

 Let $X:=\overline{C^{red}}$  and $A:=\overline{C^{red}}\setminus C^{red}$ then $X$ is a subcomplex, and $A$ is a subcomplex of $X$ by Lemma \ref{SbC}, and the result follows since 
 $C^{red}=X\setminus A$.
     \end{proof}

\begin{figure}[htbp]
\begin{center}

 \includegraphics[width=0.9\textwidth]{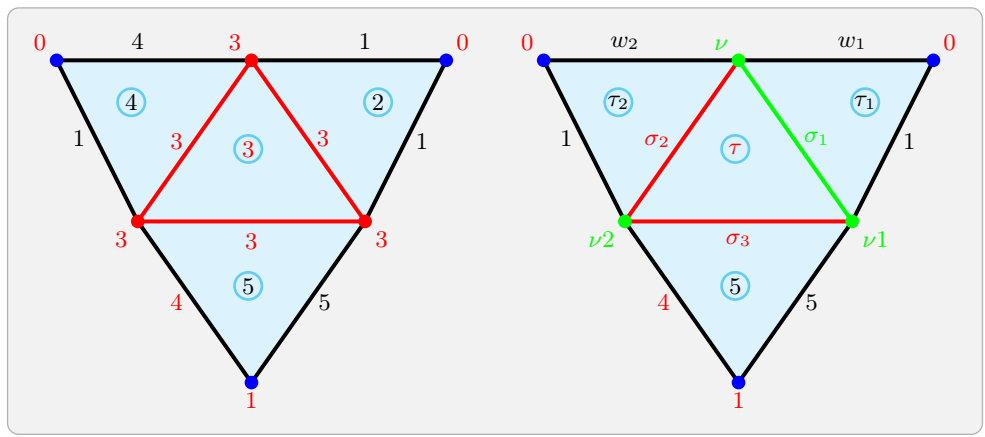}
\caption{A discrete Morse-Bott function.}\label{fig:mbt1}
\end{center}
\end{figure}

Let $b_k^{red}:=\dim \big(\ker \bo^{red}_k/\im\bo_{k+1}^{red}\big)$ then we define 
the Poincar\'e polynomial of $C^{red}$  by $$P_t(C^{red})=\sum_k b_k^{red}t^k.$$
\begin{Example}\label{expll}
Now we can change the function in Figure \ref{fig:mbt} to make it discrete Morse-Bott and  this is illustrated in Figure \ref{fig:mbt1}. Where, $C$ is the (subcomplex) collection of all the simplices with value $3$.  
 $C$ has three upward noncritical vertices and 
one upward noncritical edge, highlighted in green. Then   $\,C^{red}=\{\si_2,\si_3,\tau\}$. 
Observe that one can retrieve the Euler number of the complex just by adding that of $C^{red}$ 
to those of the critical simplices.
The contribution for $C^{red}$ is $-1$.  There are three critical 
vertices and one critical edge which all together contribute for $2$. 
Thus one has $\chi(C^{red})+2=-1+2=1=\chi(\K)$, where 
$\chi(\K)$ is the Euler number of the cell complex.
\end{Example}
Before we state our analogue of the Morse-Bott inequalities which generalizes the idea above, 
 we prove the following proposition. 

Let $n_k:=\sharp \{\si \in C^{red} | \dim \si =k\}$, then $n_k=\dim C_k(C^{red};\Z):= \dim C^{red}_k$, and let $s:=\dim C^{red}$.
\begin{Proposition}\label{MPCred}Let $C^{red}$ be a reduced collection for a discrete Morse-Bott function. Then there exists 
$r(t)$,  a polynomial in $t$ with nonnegative integer coefficients, such that 
 $$ {\sum_{k=0}}^{s}n_kt^k=P_t(C^{red}) + (1+t)r(t).$$
\end{Proposition}
\begin{proof}
If $\,\,L\,:\, V\, \rightarrow W \,$ is  a  
linear map on finitely generated $R$-modules, $R$ a principal ideal domain, then
$$\dim\, V\, =\, \dim\, \ker\,L \,+\, \dim \,\im\,L.$$

$$\dim \,C^{red}_k=n_k = \dim\, \ker \partial^{red}_k + \dim \im\partial^{red}_{k},$$  and 
$$b^{red}_k= \dim H_k\, = \,\dim \,\ker \partial^{red}_k/ \im\partial^{red}_{k+1} \,=\, \dim\, \ker \partial^{red}_k - 
\dim \,\im\partial^{red}_{k+1}$$ for $\,k= 0, 1 ,\cdots, s$, 
follows from  the  short exact sequence
$$0 \rightarrow \im \bo^{red}_{k+1} \hookrightarrow \ker \bo^{red}_k \rightarrow \dim H_k \rightarrow 0.$$
So,
\ba\nonumber
 \sum_{k=0}^s n_kt^k - \sum_{k=0}^s b^{red}_kt^k &=&\sum_{k=0}^s (\dim \ker \bo^{red}_k + \dim \im \bo^{red}_k)t^k \\\nonumber
&-& \sum_{k=0}^s(\dim \ker \bo^{red}_k- \dim \im \bo^{red}_{k+1})t^k\\\nonumber
&=&\sum_{k=0}^s (\dim \im \bo^{red}_k+ \dim \im \bo^{red}_{k+1})t^k\\\nonumber
&=& \sum_{k=0}^s(n_k-\dim \ker \bo^{red}_k)t^k +\sum_{k=0}^s(n_{k+1}-\dim \ker \bo^{red}_{k+1})t^k\\\nonumber
&=&(t+1)\sum_{k=1}^s(n_k-\dim \ker \bo^{red}_k)t^{k-1}\\\nonumber
& & \text{since}\quad n_0= \dim \ker \bo^{red}_0\,\,\text{and}\,\,n_{s+1}=0=\dim \ker \bo^{red}_{s+1}.
\ea
The proof ends by using the fact that $\,\dim \ker \bo^{red}_k \leq n_k\,$ for all $\,k=1,2,\cdots,s$.
\end{proof}
\begin{Remark}
 In general, the contribution of each noncritical pair cancels out whenever one looks at the Euler number, thus it saves a lot of time to ignore the noncritical pairs.
\end{Remark}

The idea in  Example \ref{expll}  is the content of the following analogue of the Morse-Bott inequalities.
\begin{Theorem}[Discrete Morse-Bott inequalities]\label{proeubot}
Let $f$ be a discrete Morse-Bott function on an $n$-dimensional CW complex $\K$, and let $C^{i,red}, \,\,i=1,\dots, l$ be its nonempty disjoint reduced collections that are not noncritical pairs.
Then there exists $R(t)$,  a polynomial in $t$ with nonnegative integer coefficients, such that 
\ba\label{Peoeq} 
\sum_{i=1}^l P_t(C^{i,red})=P_t(\K) + (1+t)R(t).
\ea 
\end{Theorem}
The result \rf{Peoeq} implies that, setting $t=-1$ above, 
the Euler number of each reduced collection  should always be counted with positive sign.

If we want to have the formula   
$$ \sum_{i=1}^l P_t(C^{i,red}) {t}^{ind(C^{i,red})}=P_t(\K) + (1+t)R(t),$$
then we have to put $\,\ind(C^{i,red})= 0$ for all $i=1,\cdots,l$.

\begin{Remark} 
 One of the reasons for defining the index as above is the fact that the Poincar\'e polynomial of a point in the smooth setting is given by $1$ while in the discrete setting, 
the Poincar\'e polynomial of a cell of dimension $k$ is given by $t^k$, as we shall see later on.
\end{Remark}

\begin{figure}[htbp]
\begin{minipage}[b]{0.51\linewidth}	
   \centering	
   \includegraphics[width=0.98\textwidth]{figure/mbc1.tikz}
   \caption{A discrete Morse-Bott-Conley method.}\label{fig:mbc1}
\end{minipage}\quad
\begin{minipage}[b]{0.48\linewidth}
\centering
   \includegraphics[width=0.98\textwidth]{figure/mbc2.tikz}
   \caption{A reduced collection  that is  a noncritical pair.}\label{fig:mbc2}
\end{minipage}
\end{figure}

\begin{figure}[htbp]
\begin{minipage}[b]{0.48\linewidth}	
   \centering	
   \includegraphics[width=0.98\textwidth]{figure/mbc4.tikz}
   \caption{Example of a discrete Morse-Bott-Conley method.}\label{fig:mbc4}
\end{minipage}\,\,
\begin{minipage}[b]{0.48\linewidth}
\centering
   \includegraphics[width=0.98\textwidth]{figure/mbc3.tikz}	
   \caption{Another example of a discrete Morse-Bott-Conley method.}\label{fig:mbc3}
\end{minipage}
\end{figure}

\begin{Example}
 \begin{itemize}
  \item[$1)$] In Figure \ref{fig:mbc1}, $C$ consists of the red vertex, the two green vertices, 
  the two red edges and the $2$-simplex with value $3$. Removing all the noncritical cells w.r.t. $C$ yields  that 
 $C^{red}$ consists of the two red edges and the red vertex. Thus, $P_t(C^{red})=t$;
 in addition, there are two critical vertices, thus their overall contribution is $2$.
 Hence the discrete Morse-Bott inequalities are satisfied since $\sum_i P_t(C^{i,red})=t+2$ and $P_t(\K)=1$, that is,  $R(t)\equiv 1$.
 
 \item[$2)$] In Figure \ref{fig:mbc2}, $C$ consists of the red vertex, the two red edges and the  green vertex. After 
 removing the upward or downward noncritical cells w.r.t. $C$, 
 we obtain a reduced collection that is a noncritical pair, so we do not take it 
 into consideration. Thus we only take into account the critical edge and the critical vertex, 
 and adding their contributions yield $\,t+1\,$ which is the Poincar\'e polynomial of the complex, thus, $R(t)\equiv 0$.
   \item[$3)$]In Figure \ref{fig:mbc4}, $C$ is the collection of the two red edges, 
 the green edge, the two green vertices and the $2$-simplex with value $2$. 
 $C^{red}$ consists of the two red edges and the $2$-simplex. $P_t(C^{red})=t$, 
 to this we add the contributions for the two critical vertices. We get 
 $\sum_i P_t(C^{i,red})=t+2$ and $P_t(\K)=1.$ Here, $R(t)\equiv 1$.
 \item[$3)$]In Figure \ref{fig:mbc3}, $C$ is the collection of the two red edges, 
 the green edge and the two green vertices. 
 $C^{red}$ consists of the two red edges. $P_t(C^{red})=2t$, to this we add the contributions 
 for the two critical vertices, and we get,  
 $\sum_i P_t(C^{i,red})=2t+2$, $P_t(\K)=1+t.$ Thus, $R(t)\equiv 1$.\\
  \end{itemize}
\end{Example}

The following lemma is useful for the proof of Theorem \ref{proeubot}.

\begin{Lemma}\label{lempro}
 $\sum_i \dim \ker{\bo}^{C^{i,red}}_k-\dim \ker \bo^F_k \geq 0\,$  for each $\,k\geq 1$. 
\end{Lemma}
\begin{proof}
Suppose that the reduced collections $C^{i,red}$, $i=1,\cdots,l$ are indexed such 
that $f(C^{i,red})\leq f(C^{j,red})$ for $i\leq j$. 
We know that if $\si\in \bo^F \tau$, where $\bo^F$ denotes  Forman's boundary operator, then $f(\si)<f(\tau)$. 
Using the reduced collections, 
if $\tau\in C^{i,red}$, and $\si\in \bo^F \tau$, then either $\si\in C^{i,red}$, in which case $f(\si)=f(\tau)$, or there is a 
$C^{j,red}$, s.t. $\si \in C^{j,red}$, in which case one immediately sees that $f(C^{j,red})<f(C^{i,red})$, so  $j\leq i$.\\
Also, we know that $$C_k(\K,R)=\oplus_{i=1}^l C_k(C^{i,red},R),$$ where $(R=\Z$ or $\Z_2).$ 

Let $\bo^i_k:={\bo}^{C^{i,red}}_k$, then  
$\si \in C_k(\K,R) \Rightarrow \si=\si^1+\si^2+\cdots+\si^l$ with $\si^i \in C_k(C^{i,red},R)$ where, 
\ba\nonumber
\bo^F \si^i &=&\Proj_{C^{i,red}}\bo^F \si^i + \sum_{j=1}^{i-1}\Proj_{C^{j,red}}\bo^F \si^i\\\nonumber
&=& \bo^i \si^i +\sum_{j=1}^{i-1}\Proj_{C^{j,red}}\bo^F \si^i,\quad \text{for}\quad i=1,\cdots,l,
\ea
and $\Proj$ denotes the projection map.

If $\si\neq 0$  and $\si \in \ker \bo^F$, then $\sum_i \bo^F \si^i=0.$ \\
If $\si=\si^1$ then $\bo^F \si=0 \Leftrightarrow \bo^1\si^1=0$ and we are done.\\
If at least one of the $\si^i$'$s$ is different from zero, 
 since $\Proj_{C^{j,red}}\bo^F \si^i=0$ if $j>i$, we get immediately that $\bo^l \si^l=\Proj_{C^{l,red}}\bo^F \si=0$.\\
 If $\si^l=0$, do the same for $\si^{l-1}$ and so on, until you get to $\si^1$, which we know will have to  
 be different from zero in order not to contradict the fact 
  that $\si\neq 0$. Thus at least one of the $\si^i$'$s$ must be different from zero. 
  Hence $\sum_i \dim \ker \bo^i_k\geq \dim \ker \bo^F_k$.\\\\
  Now let us suppose $\dim \ker \bo^F_k =2$.   Let   $\si_i\neq 0, $ for $i=1,2$ and $\si_1\neq \si_2$. Suppose  $\si_i \in \ker \bo^F_k$, for $i=1,2$ then 
  $\si_1=\si_1^1+\si_1^2+\cdots+ \si_1^l$ and $\si_2=\si_2^1+\si_2^2+\cdots+\si_2^l$.
  As before, we will have $\bo^l \si_1^l=0$ and $\bo^l \si_2^l=0$. 
  If $\si_1^l=\si_2^l\neq 0$, then $0\neq \bar{\si}:=\si_1+\si_2=\si_1^1-\si_2^1+\cdots+\si_1^{l-1}-\si_2^{l-1}$ satisfies  $\bo^F \bar{\si}=0$. From the previous step it follows that 
  $\sum_{i=1}^{l-1} \dim \ker \bo^i_k\geq 1.$ Using the fact that $\si_1^l \in \ker \bo^l_k$, we then get that 
  $\sum_i^l \dim \ker \bo^i_k\geq 2$.
   If  $\si_1^l=0\,\&\, \,\si_2^l \neq 0$ or  $\si_1^l\neq 0\,\&\, \,\si_2^l=0,\,\,$ we are done. 
   If $\si_1^l=0\,\&\, \,\si_2^l =0$, then do the same for $\si_1=\si_1^1+\cdots+\si_1^{l-1}$ 
   and $\si_2=\si_2^1+\cdots+\si_2^{l-1}$. 
  Hence $\sum_i^l \dim \ker \bo^i_k\geq \dim \ker \bo^F_k$.\\\\
  We assume that the statement is true for $\dim \ker \bo^F_k \leq m-1$ and we show it is true for $m$.
  Suppose $\si_i\in \ker \bo^F_k$ and are all linearly independent for $i=1,\cdots,m$, where  $\si_i=\si_i^1+\cdots+\si_i^l$.\\
  Proceeding as before, we get $\si_i^l\in \ker \bo^l_k$ for $i=1,\cdots,m$.     If $\si_i^l\neq 0$ for all $i$ and $\si_i^l \neq \si_j^l$ for all $i\neq j$ we are done. \\
  If $\si_i^l\neq 0$ for all $i$ but $\si^l_1= \si^l_2=\cdots =\si_m^l$, 
  then we get $\bar{\si}_{i}:=\si_1-\si_{i+1} \in \ker \bo^F_k$ for 
  $i=1,\cdots,m-1$. The linear independence of the $\bar{\si}_i$'$s$ follows from that of the $\si_i$'$s$, 
  and by induction hypothesis, 
  $\sum_{i=1}^{l-1} \dim \ker \bo^i_k\geq m-1$. Hence,  
  $\sum_{i=1}^{l} \dim \ker\bo^i_k\geq m=\dim \ker \bo^F_k$, also taking $\si^l_i \in \ker \bo^l_k$.\\
  If $\si_i^l\neq 0$ for all $i$ but $\si^l_1\neq \si^l_2=\cdots =\si_m^l$, 
  then we get  $\bar{\si}_{i}:=\si_2-\si_{i+2} \in \ker \bo^F_k$ for 
  $i=1,\cdots,m-2$. By induction hypothesis, $\sum_{i=1}^{l-1}\dim \ker \bo^i_k\geq m-2$. Hence,  
  $\sum_{i=1}^{l} \dim \ker\bo^i_k\geq m=\dim \ker \bo^F_k$, also adding $\si^l_1$ and $\si^l_2$ in  $\ker \bo^l_k$.\\
  The same idea is used if there are subsets $A_1,\cdots,A_s$ of $\{1,\cdots,m\}$ s.t. 
  $\si^l_i=\si^l_j$ for all $i\neq j,\,\,i,\,j\in A_p$, but $\si^l_i\neq \si^l_j$ for $\,i\in A_p$ and $\,j\in A_q$,  $\,\,p\neq q$.\\
  Indeed: $|A_1|+\cdots+|A_s|=m$, and for each $A_i$, we define $\,\bar{\si}^{A_i}_j:=\si^{A_i}_1-\si^{A_i}_{j+1}\,$ for 
  $\,j=1,\cdots,|A_i|-1.$ From the $\bar{\si}^{A_i}_j$ for $i=1,\cdots,s$ and  $j=1,\cdots,|A_i|-1$, and the induction hypothesis,
   we have: $\sum_{i=1}^{l-1}\dim \ker \bo^i_k\geq \sum_{i=1}^s(|A_i|-1)=m-s$. Taking into consideration the fact that in each $A_i$ we have 
   $\, \si^l_j \in \ker \bo^l_k$, we get $\sum_{i=1}^{l}\dim \ker \bo^i_k\geq m$. \\
     If $\si_i^l= 0$, we do the same for $\si_i^{l-1}$ and so on.
  
\end{proof}

We now have all the necessary tools to prove Theorem \ref{proeubot}.\\

\begin{proof}[\textbf{Proof of Theorem \ref{proeubot}}]
Let $f_\varepsilon\colon \K \rightarrow \R$ given by  $$f_\varepsilon (\si)=f(\si)-\frac{\varepsilon}{\dim \si +1},$$
then  $\,\,f_\varepsilon \rightarrow f\,$ as $\,\varepsilon \rightarrow 0$.\\

Let $n^i_k:=\sharp\{\si \in C^{i,red}\, | \,\dim \si=k\}$ and $s^i=\dim C^{i,red}$.\\
\\
Claim: \\
For sufficiently small $\varepsilon$, $f_\varepsilon$ is discrete Morse and\\
$\{\si^{(k)} \,\,\text{critical for}\,\, f_\varepsilon \}=\bigcup_{i=1}^l\{\si^{(k)} \in C^{i,red} \},$
that is 
\ba\label{eq(1)}
 m^{f_\varepsilon}_k = \sum_{i=1}^l n^i_k,
\ea

where 
 $m^{f_\varepsilon}_k$ is the number of critical points of $f_\varepsilon$ of dimension $k$.\\

\begin{proof}[Proof of the claim] Let $C^i$ be a collection and  $\si \in C^i$. We recall that if $\si<\tau$ is an irregular facet of $\tau$, then 
$f(\si)<f(\tau)$ and the same holds for $f_{\varepsilon}$. 
 So, it is enough to do the following at the regular facets.
 \ba\nonumber
 \sharp \{ \nu<\si | f_\varepsilon(\nu) \geq f_\varepsilon(\si)\} 
 &=& \sharp \{ \nu \in C^i, \nu<\si | f_\varepsilon (\nu) \geq f_\varepsilon (\si) \}\\\nonumber                               &+& \sharp \{ \nu \notin  C^i, \nu<\si | f_{\varepsilon}(\nu) > f_\varepsilon (\si) \}\\\nonumber
                               &=:& A_1+B_1
\ea
\ba\nonumber
 \sharp \{ \tau>\si | f_\varepsilon(\tau) \leq f_\varepsilon(\si)\} 
 &=& \sharp \{ \tau \in C^i, \tau>\si | f_\varepsilon (\tau) \leq f_\varepsilon (\si) \}\\\nonumber                               &+& \sharp \{ \tau \notin C^i, \tau>\si | f_{\varepsilon}(\tau) < f_\varepsilon (\si) \}\\\nonumber
                               &=:& A_2+B_2          
\ea
 We have, $A_1=0$ (resp. $A_2=0$) since, if  $\nu<\si$, 
 meaning that $\dim \nu=\dim \si-1$ (resp. $\tau>\si$ meaning that $\dim \tau=\dim \si+1$), 
 and both  are in the same collection that is  $f(\si)=f(\nu)$ (resp. $f(\si)=f(\tau)$), we get that 
 $f_\varepsilon(\si)>f_\varepsilon(\nu)$ (resp. $f_\varepsilon (\si)< f_\varepsilon (\tau)$).
 
 For sufficiently small $\varepsilon$, that is,  as $\varepsilon \rightarrow 0$, we have
 $B_1$ =$\sharp \{ \nu \notin  C^i, \nu<\si | f(\nu) > f (\si) \}$ and 
 $B_2$ =$\sharp \{ \tau \notin  C^i, \tau>\si | f(\tau) < f (\si) \}$,  so that, $B_1\leq1$ and $B_2\leq 1$ 
 follow from the discrete Morse-Bott condition for $f$. 
 Hence $A_1+B_1\leq 1$ and $A_2+B_2\leq 1$ which implies that  $f_{\varepsilon}$ is discrete Morse.
 
We show that 
\ba\label{eq(ii)}
\{\si^{(k)} \,\,\text{critical for}\,\, f_\varepsilon \}=\bigcup_i\{\si^{(k)} \,\,\in \,\, C^{i,red} \}.
\ea
`$\Rightarrow\,\,$' Let $\si$ be critical for $f_\varepsilon$, this implies that $B_1=0$ and $B_2=0$.
For $\varepsilon$ small enough, $B_1=0\,$ means that $\,x\,$ is not 
downward noncritical w.r.t $\,C^i$, and $\,B_2=0\,$ means that $\si$ 
is not upward noncritical w.r.t $C^i$. Therefore $\si$ should be 
in $\,C^{i,red}$. \\
`$\Leftarrow\,\,$' If  $\si\in C^{i,red}$, then $\si$ 
 was neither upward noncritical nor downward noncritical w.r.t. $C^i$ for $f$. 
 So  for sufficiently small $\varepsilon$,
 $B_1=0=B_2$
 and this implies that  $\,\si\,$ is critical for $\,f_\varepsilon$.
 \end{proof}
Now to end the proof of the  theorem: 
from Proposition
\ref{MPCred}, we have for each $~i$, 
\ba\label{eq(2)}
 {\sum_{k=0}}^{s^i} n^i_kt^k=P_t(C^{i,red}) + (1+t)r_i(t),
 \ea
where each $r_i(t)$ is a polynomial in $t$ with nonnegative integer coefficients.\\
The function $f_\varepsilon$ is discrete Morse on $\K$ for sufficiently small $\varepsilon$, so from the  Morse inequalities we have:
\ba\nonumber
& & \sum_{k=0}^n m^{f_\varepsilon}_k t^k= \sum_{k=0}^n b_k t^k + (1+t)r^1(t)\\
\nonumber
&\Rightarrow & \sum_{k=0}^n (\sum_{i=1}^l n^i_k )t^k= \sum_{k=0}^n b_k t^k + (1+t)r^1(t) \quad \text{from}\,\,\rf{eq(1)}\\
\nonumber
& \Rightarrow& \,\, \sum_{i=1}^l P_t(C^{i,red})+ (1+t)r_i(t)
= \sum_{k=0}^n b_k t^k + (1+t)r^1(t) \quad \text{from}\,\,\rf{eq(2)}\\
\nonumber
& \Rightarrow &\sum_{i=1}^l P_t(C^{i,red}) =P_t(\K)+ (1+t)R(t).
\ea
Note that $$\,\,r^1(t)=\sum_{k=1}^n (m^h_k-\dim \ker \bo^F_k)t^{k-1}\,\,\text{ and}\,\, 
r_i(t)={\sum}^{s^i}_{k=1}(n_k^i- \dim \ker \bo_k^{C^{i,red}})t^{k-1},$$  where
$\,\,\bo_k^{C^{i,red}}:={\bo^c_k}_{| C^{i,red}},\,\,$. Thus 
$R(t)= \sum_{k=1} (\sum_i \dim \ker\bo^{C^{i,red}}_k-\dim \ker \bo^F_k)t^{k-1}$,   and the result follows from Lemma \ref{lempro}. 
\end{proof}

After realizing that an approach of Floer's in the case of a discrete Morse function is possible 
one would like to ask the question 
if something similar can be done  
using a discrete Morse-Bott function. \\

\textbf{Question:} \textit{is it possible to do some kind of Floer-related theory on a complex having a discrete Morse-Bott function?}\\

It turns out that, 
 we cannot consider the whole reduced collection as one critical object. Thus the only solution is: 

\textbf{Solution:} \textit{we approximate the discrete Morse-Bott function on each collection to get a discrete Morse function; this is always possible as seen before, so that   
any upward or downward noncritical cell w.r.t. a given collection will not be  critical after the approximation. 
In this way, all the cells in the reduced collections will be critical just by making the approximating function to be a trivial discrete Morse function.}\\

The above solution shows that  the only way we can have a boundary operator on some complex on which 
is defined our discrete Morse-Bott function is by approximating 
to get a discrete Morse function. We now ask:  having a discrete Morse-Bott function, is it possible to develop  any other theory
 using a more  dynamics-related method 
to get some insight on the Euler number?  This is what we will be investigating in the last part of this paper.

Now we use the discrete Morse-Bott function to derive a similar notion to smooth Conley theory, this is our version of discrete Conley theory.

\section{Discrete Morse-Bott-Conley theory}\label{sec5.4}

We recall, see  \cite{J2}, \cite{KMi} and \cite{CC}, that Conley theory  uses 
the idea of exit sets and isolating neighborhoods of isolated invariant sets to derive the Euler
characteristic of the manifold at hand. For each isolated invariant set, the pair consisting of its isolating neighborhood and exit set is called its index pair. 
 The Euler number of the manifold is then obtained by summing for all isolating invariant sets, the alternating sums of the dimensions 
of the homologies of  the index pairs.

The scope of Conley theory is more general than that of Morse theory,  because 
 if we consider the discrete vector field originating from some discrete Morse 
function, then this combinatorial vector field does not admit any closed orbit. This makes things  easy in the sense that, the only critical objects under consideration are the 
critical cells of the discrete Morse function, and we will show below that for a critical cell of dimension $k$, 
the homological Conley index is just given by $t^k\,$, or $\,(-1)^k$ for $t=-1$.

Forman, in \cite{For3}, made things more interesting by considering a combinatorial vector field in general 
(it need not originate from some discrete Morse function and can therefore admit closed orbits).
In general, a combinatorial vector field on a CW complex yields a disjoint collection  $C_r$ of rest points or closed orbits,
 where the rest points are exactly the critical cells. For each $C_r$, let $\overline{C_r}$ 
be the union of all the cells in $C_r$ with the ones in their boundaries, and let $\widetilde{C_r} \,:=\, \overline{C_r}\setminus C_r$. 
The pair $(\overline{C_r},\widetilde{C_r})$ is taken to be  an index pair  for the corresponding  $C_r$. Forman showed that with
$$m_i\,:=\,\sum_{C_r} \dim H^i(\overline{C_r},\widetilde{C_r};\mathbb Z),$$ the Morse inequalities are satisfied.

Now we want to go beyond a combinatorial vector field and we do this by considering 
 our discrete Morse-Bott function on a CW complex $\K$. The vector field originating from this 
function inside each collection need not be a combinatorial vector field 
as a cell  in the collection can have more than one incoming and/or outgoing arrow,  
but between the collections there should not be any closed orbit.

We shall consider as isolated invariant sets the reduced collections (excluding the noncritical pairs), and their isolating neighborhood will be the subcomplex generated by the 
isolated invariant set. The exit set will be 
the part of $\bo^{top}N$ where the values of the function are smaller than on the isolated invariant set. This is the content of the following definition.

\begin{Definition}
\begin{itemize}
 \item [$1)$] $I$ is said to be an \textbf{isolated invariant} set if it is a reduced collection that is not a noncritical pair.
  \item[$2)$] An \textbf{isolating neighborhood} $N$ for $I$ is the union of all the cells in $I$ together with all the cells in their boundaries. That is,  
  $$N=\bigcup_{\si \in I}\bar{\si}.$$
  \item[$3)$] The \textbf{exit set} for the flow from $I$ is just given by  
 the part of $N$ not in $I$ where the values of $f$ are smaller than
 or equal to the value on $I$. In other words 
$$E=\{\si\in N\setminus I\,\,| \,\,f(\si)\leq f(\tau) \,\,\text{for}\,\, \tau \in I\}.$$
\item[$4)$] We call $(N,E)$ an index pair for $\,I$.
\end{itemize}
\begin{Remark}
\begin{itemize}
\item[$1)$] The reduced collections $I^{red}_i$ are such that there is no path (following the arrows) that moves from a cell in $I^{red}_i$ to another cell 
 outside of $I^{red}_i$, meaning that each of them  is  invariant. 
 Since they constitute the building block for computing the Poincar\'e polynomial of the CW complex, 
 we consider them  to be isolated. Hence each collection $I^{red}_i$ is an isolated invariant set.
 \item[$2)$] In the definition of the exit set above, taking $I=C^{red}$, the cells in $E$  that have the same value with the ones in $I$ are exactly 
 those cells that are upward  noncritical w.r.t. $C$. Indeed, we know from Lemma \ref{SbCb} that  
any $\si\in N\setminus I$ satisfying $f(\si)=f(\tau)$ for $\tau \in I$ cannot be downward noncritical. Thus,
\ba\nonumber
E = \Bigg\{  \begin{aligned}\si\in N\setminus I \quad \text{s.t.} \,\,\, 
\text{ either}\,\, f(\si)<f(\tau)\,\, \text{for}\,\, \tau\in I,\,\,\,\text{ or}\\ 
    f(\si)=f(\tau)\,\,\text{ and}\,\,\, \si\, \,\text{is upward noncritical w.r.t.}\,\,C \end{aligned}
    \Bigg\}.
\ea  
Later on, see the proof of Theorem \ref{TheoBC}, we will show that $E=N\setminus I$.
\end{itemize}
\end{Remark}

\end{Definition}

\begin{Definition}
The \textbf{topological Conley index} of an isolated invariant set $I$ 
is the homotopy type of $N/E$  and the \textbf{homological Conley} index of $I$ is   the polynomial
$$C_t(I)\,:=\,\sum_k \dim H_k(N,E;\mathbb Z)t^k.$$
\end{Definition}

We then have the following lemma.

\begin{Lemma}\label{Lem305}
 If  $I=C^{red}\,$ with $\,\,\sharp C^{red}=1$ of dimension $k$, then its homological Conley index is    $t^k$.
\end{Lemma}

\begin{proof}

Let $\,I=\si^{(k)}$,  then 
$\,\,N=\bar{\si}$, and  
 $\,\,E= \bo^{top} \bar{\si} $,  
 where $\bo^{top}$ denotes the topological boundary. 
 Thus we obtain an index pair $(\bar{\si},\bo^{top}(\si))\,$ for $\sigma$, and the homotopy type of this index pair 
is just  that of a $k$-dimensional 
closed disc with its boundary identified, which is that of the pair $(S^k,pt)$. Thus
the Conley index is given by $$C_{t}(\sigma^{(k)})\, = \sum_i\,dim \,H_i(\Se^k,pt)t^i= t^k,$$ since $H_i(\Se^k,pt)=1$ for $i=k$ and $H_i(\Se^k,pt)=0$ for $i\neq k$.
\end{proof}

In general, if $I=C^{red}$ is such that $C^{red}\neq C$, then there are elements that are upward noncritical w.r.t. $C$,  
and (in order not to contradict the definition of a discrete Morse-Bott function) there is a one to one correspondence between those elements 
and those out of $C$ making them upward noncritical, this is the content of Lemma \ref{Lem02}. Thus we have disjoint noncritical pairs 
$\{\si<\tau \}$,  $\si \in C$ and  $w(\si)=\tau \notin C$ with $f(\si)\geq f(\tau)$. 
We call such a $\si$ an \textbf{exit cell}.\\

\begin{Remark}
 If $I=C^{red}$, then the exit set of 
$I$ denoted $E$ is the union of all the cells in $N$ whose values are smaller than the value in $C^{red}$, together with  all the 
 exit cells of $C$ that are contained in $N$.
\end{Remark}

 We have  the following analogue of the result in the smooth setting:

\begin{figure}[htbp]
\begin{center}

 \includegraphics[width=0.98\textwidth]{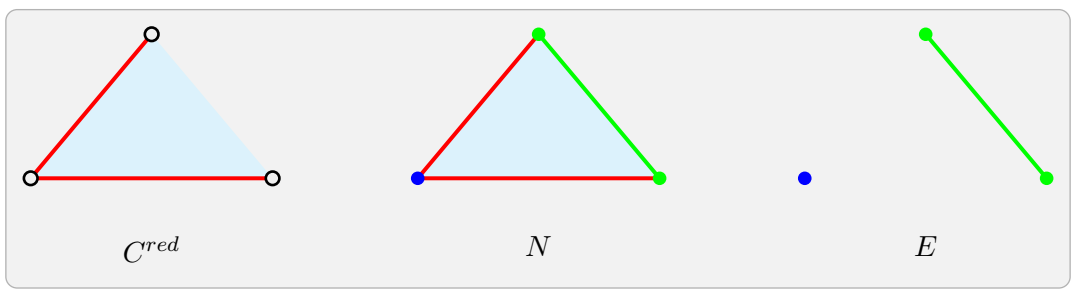}
\caption{A reduced collection and its index pair.}\label{fig:mbc33}
\end{center}
\end{figure}

\begin{figure}[htbp]
\begin{center}

 \includegraphics[width=0.98\textwidth]{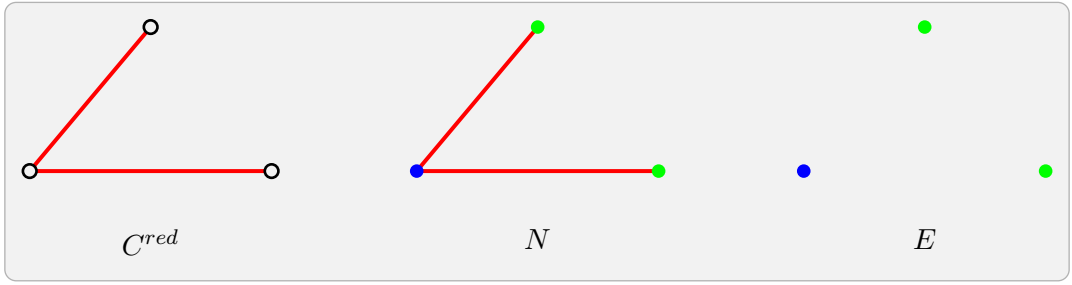}
\caption{Another reduced collection and its index pair.}\label{fig:mbc44}
\end{center}
\end{figure}
\begin{Theorem}\label{TheoBC}Let $f$ be a discrete Morse-Bott function having isolated invariant sets 
$I_1,\cdots,I_l$, then  there exists  $R(t)$, a polynomial in $t$ with nonnegative integer coefficients, such that \ba\label{coneu}
P_t(\K) +(1+t)R(t)\,=\,\sum_{j=1}^l\, C_{t}(I_j). \ea
\end{Theorem}

Recall that $P_t(N,E)=P_t(N/E)-1$.

\begin{Example}
 \begin{itemize}
  \item [$1)$] For the discrete Morse-Bott function in Figure \ref{fig:mbc4}, we get the respective index pair in Figure \ref{fig:mbc33}. Thus $P_t(N,E)=t$. 
Observe that this can be achieved geometrically by performing 
  $N/E$ and using the fact that $P_t(N,E)=P_t(N/E)-1.$ Hence, $P_t(\K)=1$ and $\sum_jC_t(I_j)= t + 2$, this 
  implies that $R(t)\equiv 1$.
  \item[$2)$]Using Figure \ref{fig:mbc3}, the respective index pair for the discrete Morse-Bott function is given by Figure \ref{fig:mbc44}. In this case, 
 we get $P_t(N,E)=2t$. Thus  $P_t(\K)=1+t$ and $\sum_jC_t(I_j)= 2t + 2$, this 
  implies that $R(t)\equiv 1$.
 \end{itemize}

\end{Example}

\begin{proof}[\textbf{Proof of Theorem \ref{TheoBC}}]
\begin{itemize}
\item We have shown above, see Lemma \ref{Lem305}, that whenever $I$ is a singleton, a cell of dimension 
$k$, $$C_{t}(I) \,= \,t^k.$$ 
\item If
$I=C^{red}$ where $\sharp C^{red}>1$,  we only need to show that $C^{red}=N\setminus E$ (which holds for all $C^{red}$). \\
The fact that $C^{red} \subseteq N\setminus E$ follows from the definitions of those sets.\\
We show that $N\setminus E \subseteq C^{red} .$\\
 Let $\si \in \,N\setminus E$, then  $\si$ is either in $C^{red}$ or it is a face of an element in $C^{red}$. If $\si \in C^{red}$, we are done.
If  there is a $\tau \in C^{red}$ s.t. 
$\si<\tau$, then from the definition of $C^{red}$ and the fact that  $\si$ is not in $E$, we must have $f(\si)=f(\tau)$ and $\si$ is not upward noncritical.
Lemma \ref{SbCb} also establishes the fact that $\si$ cannot be downward noncritical.
 Thus, $f(\si)=f(\tau)$ and $\si$ is neither upward nor downward noncritical, this implies $\si$ is in $C^{red}$.\\
Thus $N:=N(I)$ is a subcomplex by definition, and $E:=E(I)=N\setminus I:=\overline{C^{red}}\setminus C^{red}$ is a  subcomplex by Lemma \ref{SbC}. 
Hence, the relative homology is well defined. 
Thus, 
$$C_{t}(I)= \sum_{k}\dim H_k(N(I),E(I);\Z)t^k=P_t(I),$$
and the  second equality follows from   Theorem \ref{proeubot}  since $C^{red}=N\setminus E$ and both $N:=N(I_j)$ and $E:=E(I_j)$ are subcomplexes, indeed
\ba\nonumber
\sum_{j=1}^l\, C_{t}(I_j)&=&\sum_{j=1}^l\sum_{k}\dim H_k(N(I_j),E(I_j);\Z)t^k \\\nonumber
&=& \sum_{j=1}^lP_t(I_j)=P_t(\K) +(1+t)R(t).\ea

\end{itemize}
\end{proof}
A consequence of the proof above is the following:
\begin{Lemma}
If $(N,E)$ is an index pair for $C^{i,red}$, then 
 $$\chi(C^{i,red})=\chi(N,E)=\chi(N)-\chi(E).$$
 \end{Lemma}
\begin{Remark}
 The lemma above is not true in general if we consider the Poincar\'e polynomials instead of the Euler numbers.
\end{Remark}

\begin{figure}[H]
\begin{center}

 \includegraphics[width=0.92\textwidth]{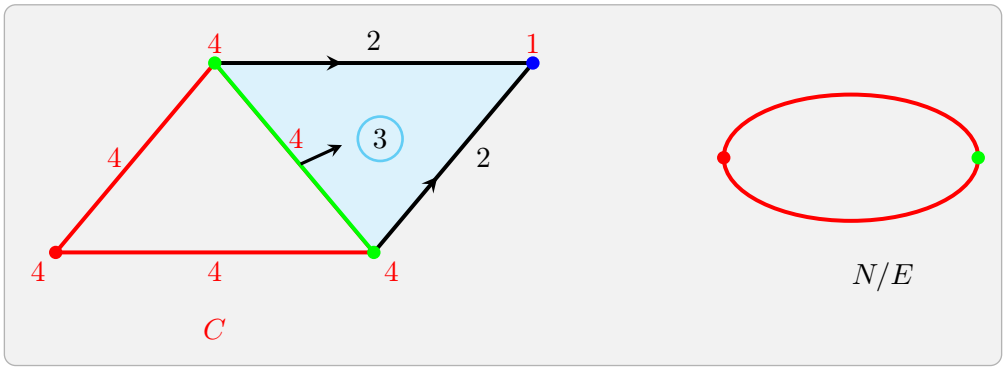}
\caption{The corresponding quotient space to an index pair.}\label{fig:mc2}
\end{center}
\end{figure}

\begin{figure}[H]
\begin{center}

 \includegraphics[width=0.92\textwidth]{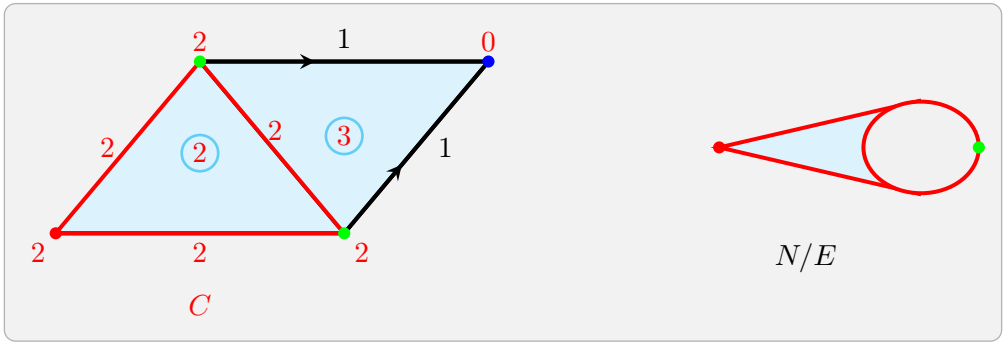}
\caption{The quotient space corresponding to an index pair.}\label{fig:mc1}
\end{center}
\end{figure}
\begin{Example}
 \begin{itemize}
  \item [$1)$] Using Figure \ref{fig:mc2}, $C^{red}$ consists of the two red edges and the red vertex, $N$ is the subcomplex generated by the two red edges, 
  so that $E$ is the collection of the two green vertices. $P_t(N,E)=P_t(N/E)-1=t$. The critical vertex contributes $1$. Thus 
  $\sum _j C_t(I_j)= t+ 1=P_t(\K)$ that is $R(t)\equiv 0$.
  \item[$2)$] In Figure \ref{fig:mc1}, $C^{red}$ takes all the elements of $C$ except the two green vertices. 
  $N=C$ and $E$ is the two green vertices. $P_t(N,E)=t$, there is one critical $2$-simplex whose contribution is $t^2$, and one critical vertex whose contribution is $1$. 
  Hence $\sum_jC_t(I_j)=t^2+t+1$, but $P_t(\K)=1$ this implies that $R(t)=t$.
 \end{itemize}

\end{Example}

\noindent\textbf{\Large Acknowledgement}\\
The author thanks J. Jost, J. Portegies, R. Matveev, D. Tran,    and R. Wu for useful remarks, suggestions and/or discussions.\\
This work was supported by a stipend
from the International Max Planck Research School, IMPRS.

\end{document}